\documentclass{SCDG2019}

\usepackage{euscript}

\newcommand{\beq}{\begin{equation}}
\newcommand{\beqnt}{\begin{equation}\notag}
\newcommand{\eeq}{\end{equation}}
\newcommand{\bpr}{\begin{proof}}
\newcommand{\epr}{\end{proof}}
\newcommand{\fref}[1]{{\rm(\ref{#1})}}                      %%%%%%%%%% formula ref      %%%%%
\newcommand{\RA}{\ensuremath{\mathbb R}}%  действительные числа
\newcommand{\mydef}{\triangleq}% равно по определению
\renewcommand{\ge}{\geqslant}% Больше равно
\renewcommand{\le}{\leqslant}% Меньше равно
\newcommand{\myles}{\mathrel\prec}% знак строгого частичного порядка
\newcommand{\myand}{\ensuremath{\&}}% и
\newcommand{\myimp}{\ensuremath{\Rightarrow}}% имп
\newcommand{\myll}{\ensuremath{\forall}}% квантор всеобщности
\newcommand{\myxst}{\ensuremath{\exists}}% квантор существования
\newcommand{\icA}{\ensuremath{\operatorname{\mathbf A}}}% ОПП
\newcommand{\icI}{\ensuremath{\mathbf{I}}}% время
\newcommand{\icX}{\ensuremath{\mathbf{X}}}% прос-во
\newcommand{\icP}{\EuScript{P}}% булеан 
\newcommand{\icPp}{\EuScript{P^\prime}}% булеан без пустого множества
\newcommand{\icPi}[3]{\ensuremath{\Pi(#1\,|\,#2,#3)}}% П
\newcommand{\icDc}{\ensuremath{\mathcal D}}% Dкалиграф
\newcommand{\icDb}{\ensuremath{\mathbf D}}% жирн.D
\newcommand{\icN}{\ensuremath{\EuScript N}}% фаз.огран
\newcommand{\icC}{\ensuremath{\mathbf C}}% траектор.
\newcommand{\icY}{\ensuremath{\mathbf{Y}}}% прос-во неопред.
\newcommand{\icOm}{\ensuremath{\Omega}}% неопред.
\newcommand{\icS}{\ensuremath{\EuScript S}}% система
\newcommand{\icMa}{\ensuremath{\mathbb M}}% к/с
\newcommand{\naZ}{\ensuremath{\mathrm{Z}}}% реал.ответов
\newcommand{\naZo}[2]{\naZ_0\setbox0=\hbox{\text{$#1$}}\ifdim\wd0=0pt\else\res{#1}{#2}\fi}%  росток ответов
\newcommand{\naOmo}[2]{\icOm_0\setbox0=\hbox{\text{$#1$}}\ifdim\wd0=0pt\else\res{#1}{#2}\fi}%  росток ответов
\newcommand{\res}[2]{\ensuremath{(#1\,|\,#2)}}% сужение #1 на #2
\newcommand{\sres}[2]{\ensuremath{\text{\large\textbf{\textup (}}#1\,|\,#2\text{\large\textbf{\textup )}}}}% сужение множества #1 на #2
\newcommand{\sysp}{{\bf SP}_{(t_0,x_0)}}

\newcommand{\eissue}{{\noindent\tiny Compiled {\tt \number\day.\number\month.\number\year} from {\tt \jobname.tex}}\par}
\newcommand{\rissue}{{\noindent\tiny Получено \today\ из файла {\tt \jobname.tex}}\par}

\begin{document}
\UDK{517.977} % УДК, соответствующий тематике вашей статьи

\Title{О динамической задаче удержания с неразложимым множеством помех}
{Д.А.~Серков} %% Список авторов для аключения в содержание сборника

\Auth{Д.А.~Серков\org1}

% Места работы и адреса авторов
% Рекомендуем использовать официальное сокращенное название организации, например, ИПУ РАН, ВолГУ, ИММ УрО РАН и т.п.
\Addr{\org1 Институт математики и механики им. Н.Н.Красовского УрО РАН, Екатеринбург, Россия, serkov@imm.uran.ru}

\Abstract{Для абстрактной динамической системы рассмотрена задача удержания движений в заданном подмножестве пространства историй движения.
Изучается случай неразложимого множества помех. 
Предложена конструкция разрешающей квазистратегии, опирающаяся на метод программных итераций.
}
\Keyword{задача удержания, неразложимое множество помех, квазистратегия.}

% При необходимости допускается деление текста на разделы с помощью команды \section{}.

%%%%%%%%%%%%%%%%%%%%%%%%%%%%%%%%%%%%%%%%%%%%%%%%%%%%%%%%%%%
\section*{Введение}
Рассмотрена задача удержания движений абстрактной динамической системы в заданном множестве --- простой случай позиционной дифференциальной игры \cite{KraSub74}.
Решение задачи ищется во множестве квазистратегий.
Новым в постановке является отказ от свойства разложимости \cite{Rock-PJM-1968} множества помех. 
В задачах управления это свойство определяется как возможность <<склейки>> двух допустимых помех в любой момент времени. 
Пример отсутствия разложимости множества помех дает случай непрерывных (постоянных) помех.
В качестве основы предлагаемого решения использован непрямой метод программных итераций \cite{Chentsov76DAN227}.
Отсутствие топологических требований компенсируется увеличением количества итераций \cite{Ser_TRIMM2017}.

%%%%%%%%%%%%%%%%%%%%%%%%%%%%%%%%%%%%%%%%%%%%%%%%%%%%
\section{Постановка задачи}

\paragraph{Обозначения.}
Через $\icP(T)$ (через $\icP'(T)$) обозначим все (все непустые) подмножества множества $T$.
Если $A$ и $B$ --- непустые множества, то $B^A$ есть множество всех отображений из множества $A$ в множество $B$. 
Если при этом $f\in B^A$ и $C\in\icP'(A)$, то $\res{f}{C}\in B^C$ есть сужение $f$ на множество $C$: $\res{f}{C}(x)\mydef f(x)$ $\myll x\in C$. 
В случае, когда $F\in\icP'(B^A)$, полагаем $\sres{F}{C}\mydef\{\res{f}{C}:f\in F\}$.

\paragraph{Динамическая система.}
Фиксируем непустое подмножество \icI\ вещественной прямой \RA\ в качестве аналога временного интервала и непустое множество \icX, задающее область изменения пространственных переменных. 
Если $t\in\icI$, то $\icI^t\mydef\{\xi\in\icI\mid\xi\le t\}$ и $\icI_t\mydef\{\xi\in\icI\mid\xi\ge t\}$.
Выделяем множество $\icC\in\icP'(\icX^\icI)$ в качестве аналога траекторий. 
Фиксируем непустые множества \icY\ и $\icOm\in\icP'(\icY^\icI)$ --- допустимые реализации помехи. 
Полагаем $\icDc\mydef\icI\times\icC\times\icOm$ --- фазовое пространство системы, понимаемое как пространство состояний управляемого процесса.
Для любых $t\in\icI$, $x\in\icC$ обозначим $\naZo x{\icI^t}\mydef\{x'\in\icC\mid\res{x'}{\icI^t}=\res{x}{\icI^t}\}$.

Фиксируем (в качестве аналога <<системы>>) отображение $\icS:\icDc\mapsto\icPp(\icC)$ такое, что $\myll t\in\icI$, $\myll\tau\in\icI_t$ $\myll x,x^\prime\in\icC$ и $\myll\omega,\omega'\in\icOm$
\beq\label{sys1}
\icS(t,x,\omega)\in\icPp(\naZo{x}{\icI^t}),
\eeq
\beq\label{sys2}
(\res{x}{\icI^t}=\res{x^\prime}{\icI^t})\myimp\left(\icS(t,x,\omega)=\icS(t,x^\prime,\omega)\right),
\eeq
\beq\label{sys3}
(h\in\icS(t,x,\omega))\myimp(h\in\icS(\tau,h,\omega)),
\eeq
\begin{multline}\label{sys4}
\Bigl(\sres{\icS(t,x,\omega)}{\icI^\tau}=\sres{\icS(t,x,\omega^\prime)}{\icI^\tau}\\
\myand(h\in\icS(t,x,\omega))\myand(h^\prime\in\icS(\tau,h,\omega^\prime))\Bigr)\\
\myimp(h^\prime\in\icS(t,x,\omega^\prime)).
\end{multline}
Если $(t,x,\omega)\in\icDc$, то $\icS(t,x,\omega)$ есть множество всех траекторий <<системы>> \fref{sys1}--\fref{sys4}, отвечающих начальной истории $x$ до момента $t$ и действию помехи $\omega$ после момента $t$.

Фиксируем начальную историю $(t_0,x_0)\in\icI\times\icC$.
Все дальнейшие построения проводятся с целью формулировки и решения задачи удержания для этой начальной истории.
Выделим в $\icDc$ множество $\sysp$ всех состояний управляемого процесса, возникших в системе из начальной истории $(t_0,x_0)$ при реализации всех возможных помех: $\sysp\mydef\{(t, x, \omega)\in\icDc\mid t\in\icI_{t_0}\ \res{x}{\icI^t}\in\sres{\icS(t_0,x_0,\omega)}{\icI^t}\}$.
Для каждого $(t,x,\omega)\in\sysp$ определим множество допустимых (совместимых) помех: $\icOm(t,x,\omega)\mydef\{\omega'\in\icOm\mid \sres{\icS(t_0,x_0,\omega)}{\icI^t}=\sres{\icS(t_0,x_0,\omega')}{\icI^t}\}$.
Тогда  $\icOm(t_0,x_0,\omega)=\icOm$ для всех $\omega\in\icOm$.

%%%%%%%%%%%%%%%%%%%%%%%%%%%%%%%%%%%%%%%%%%%%%%%%%%%%%%%%%%%%%%%%%
\paragraph{Процедуры управления и задача удержания.}

Управляющая сторона для формирования пучка траекторий использует непустозначные неупреждающие мультифункции из $\icP(\icC)^\icOm$.
Итак, для $(t,x,\omega)\in\sysp$ определим множество $\icMa_{(t,x,\omega)}$ допустимых процедур управления --- квазистратегий:
\begin{multline*}
\icMa_{(t,x,\omega)}\mydef\Big\{\alpha\in\!\!\!\!\!\!\prod_{\omega'\in\icOm(t,x,\omega)}\!\!\!\!\!\!\icPp(\icS(t,x,\omega'))\ \mid\ \myll\omega_1,\omega_2\in\icOm(t,x,\omega)\ \myll\tau\in\icI\\\
\Big(\sres{\icS(t_0,x_0,\omega_1)}{\icI^\tau}=\sres{\icS(t_0,x_0,\omega_2)}{\icI^\tau}\Big)\\
\myimp\Big(\sres{\alpha(\omega_1)}{\icI^\tau}=\sres{\alpha(\omega_2)}{\icI^\tau}\Big)\Big\}.%\label{quasy}
\end{multline*}

Пусть множество $\icDb\in\icPp(\icI\times\icC)$, описывающее заданные фазовые ограничения, удовлетворяет условиям: $(t_0,x_0)\in\icDb$, $(t,x)\in\icDb)\myimp(\{t\}\times\naZo x{\icI^t}\subset\icDb)$.
Определим на его основе множество $\icN$ вида $\icN\mydef(\icDb\times\icOm)\cap\sysp$ и поставим целью управления удержание состояний системы во множестве $\icN$.
Мы будем считать ее достижимой для начального состояния $(t,x,\omega)$, если для некоторой квазистратегии $\alpha\in\icMa_{(t,x,\omega)}$ выполнены включения $(\tau,h,\nu)\in\icN$, $\myll\tau\in\icI_t$, $\myll h\in\alpha(\nu)$, $\myll\nu\in\icOm(t,x,\omega)$.

Для начальной позиции $(t_0,x_0)$ это эквивалентно удержанию текущих позиций системы в пределах $\icDb$ при любых реализациях помехи $\nu\in\icOm$.

%%%%%%%%%%%%%%%%%%%%%%%%%%%%%%%%%%%%%%%%%%%%%%%%%%%%%%%%%%%%%%%
\section{Квазистратегии, разрешающие задачу удержания.}

%%%%%%%%%%%%%%%%%%%%%%%%%%%%%%%%%%%%%%%%%%%%%%%%%%%%
\paragraph{Оператор программного поглощения и его итерации.}

Пусть $H\in\icP(\sysp)$, $(t,x,\omega)\in\sysp$ и $\nu\in\icOm(t,x,\omega)$. 
Обозначим
\beq\label{Pi}
\icPi{\nu}{(t,x,\omega)}{H}\mydef\{h\in\icS(t,x,\nu)\mid(\tau,h,\nu)\in H\ \myll \tau\in\icI_t\}.
\eeq
В терминах \fref{Pi} введем оператор $\icA\in\icP(\sysp)^{\icP(\sysp)}$ (программного поглощения): для любого $H\in\icP(\sysp)$
$$%\beq\label{OPP}
\icA(H)\mydef\{(t,x,\omega)\in H\mid\icPi{\nu}{(t,x,\omega)}{H}\neq\varnothing\ \myll\nu\in\icOm(t,x,\omega)\}.
$$%\eeq

Для произвольного ординала $\alpha$, следуя методу трансфинитной индукции \cite[п. I.3]{Engelking1986}, введем $\alpha$-итерацию $\icA^\alpha\in\icP(\sysp)^{\icP(\sysp)}$ оператора $\icA$:
при $\alpha=0$ положим $\icA^0(H)\mydef H$, $\myll H\in\icP(\sysp)$;
если $\alpha$ имеет предшественника (пусть это ординал $\gamma$), примем $\icA^\alpha\mydef\icA\circ\icA^\gamma$;
для предельного ординала $\alpha$, положим $\icA^\alpha(H)\mydef\cap_{\beta\myles\alpha}\icA^\beta(H)$, $\myll H\in\icP(\sysp)$.

Обратимся к вопросу о разрешимости в выбранном классе квазистратегий задачи удержания. 
Пусть ординал $\sigma$ строго больше мощности множества $\icN$. 
Тогда справедливы следующие утверждения.

%%%%%%%%%%%%%%%%%%%%%%%%%%%%%%%%%%%%%%%%%%%%%%%%%%%%
\begin{stat}\label{ass-P-in-M}
Для любой $(t,x,\omega)\in\icA^\sigma(\icN)$ выполняется включение $\icPi{\cdot}{(t,x,\omega)}{\icA^\sigma(\icN)}\in\icMa_{(t,x,\omega)}$.
\end{stat}

%%%%%%%%%%%%%%%%%%%%%%%%%%%%%%%%%%%%%%%%%%%%%%%%%%%%
\begin{theo}\label{ass-MzhN-in-Ai}
Выполняется равенство
\begin{multline*}%\label{MzhN=Ai}
\icA^\sigma(\icN)=\{(t,x,\omega)\in\icN\mid\myxst\alpha\in\icMa_{(t,x,\omega)}:\\
 (\tau,h,\nu)\in\icN\ \myll \tau\in\icI_t\ \myll h\in\alpha(\nu)\ \myll\nu\in\icOm(t,x,\omega)\}.
\end{multline*}
\end{theo}
%%%%%%%%%%%%%%%%%%%%%%%%%%%%%%%%%%%%%%%%%%%%%%%%%%%%
В силу теоремы \ref{ass-MzhN-in-Ai} исходная задача удержания разрешима, если и только если $(t_0,x_0,\omega_0)\in\icA^\sigma(\icN)$ для некоторого $\omega_0\in\icOm$;
при этом (см. утверждение \ref{ass-P-in-M}), в случае ее разрешимости квазистратегия $\icPi{\cdot}{(t_0,x_0,\omega_0)}{\icA^\sigma(\icN)}$ реализует это решение.

\end{document}